\documentclass[11pt,a4paper]{amsart}
\usepackage{amsmath,amsthm}
\newcommand{\bo}{\boldsymbol}
\DeclareMathOperator{\sgn}{sgn}
\newcommand{\zz}{\mathbb Z}
\newcommand{\nn}{\mathbb N}
\newcommand{\zzh}{\,\backslash\!\!\!\!\mathbb Z}
\swapnumbers
\theoremstyle{definition}
\newtheorem{point}{}[section]
\theoremstyle{plain}
\newtheorem{lemma}[point]{Lemma}
\newtheorem{theorem}[point]{Theorem}

\newtheorem{prop}[point]{Proposition}
\newcommand{\marginextend}[1]{ \addtolength{\oddsidemargin}{-#1}  \addtolength{\evensidemargin}{-#1}
  \addtolength{\textwidth}{#1}\addtolength{\textwidth}{#1}}
\newcommand{\updownextend}[1]{ \addtolength{\topmargin}{-#1}  \addtolength{\textheight}{#1}
\addtolength{\textheight}{#1}}
\marginextend{.9cm}
\updownextend{0.0cm}
\author{Gyula Lakos}
\title{Factorization of Laurent series over commutative rings}
\address{Alfr\'ed R\'enyi Institute of Mathematics, Hungarian Academy of Sciences,
Budapest, Pf. 127, H--1364, Hungary}
\address{Institute of Mathematics, E\"otv\"os University, P\'azm\'any P\'eter s.~1/C,
Budapest, H--1117, Hungary}
\email{lakos@cs.elte.hu}
\begin{document}
\begin{abstract}
We generalize the Wiener-Hopf factorization of Laurent series to
more general commutative coefficient rings,
and we give explicit formulas for the decomposition.
We emphasize the algebraic nature of this factorization.
\end{abstract}
\maketitle
\section*{Introduction}
Suppose that we have a Laurent  series $a(z)$ with rapidly decreasing coefficients
and it is multiplicatively invertible, so its multiplicative inverse $a(z)^{-1}$ allows a Laurent series expansion
also with rapidly decreasing  coefficients; say
\[a(z)=\sum_{n\in\zz}a_nz^n,\qquad\text{and}\qquad
a(z)^{-1}=\sum_{n\in\zz}b_nz^n.\]

Then it is an easy exercise in complex function theory that there exists a factorization
\[a(z)=a^-(z) \tilde a(z) a^+(z)\]
such that
\[a^-(z)=\sum_{n\in\nn} a_n^- z^{-n},\qquad\quad \tilde a(z)=\tilde a_p z^p,\quad\qquad
a^+(z)=\sum_{n\in\nn} a_n^+ z^{n},\]
\[a^-(z)^{-1}=\sum_{n\in\nn} b_n^- z^{-n},\qquad \tilde a(z)^{-1}=\tilde a_p^{-1} z^{-p},\qquad
a^+(z)^{-1}=\sum_{n\in\nn} b_n^+ z^{n},\]
\[a_0^-=1,\qquad  b_0^-=1,\qquad p\in \zz,\qquad a^+_0=1,\qquad  b^+_0=1,\]
where $a_n^+$, $a_n^-$, $b_n^+$, $b_n^-$ are rapidly decreasing.

If $G$ denotes the unit group of Laurent series with rapidly decreasing coefficients
then we find that there is an inner direct product decomposition $G=G^-\times \tilde G\times G^+$
to strictly  antiholomorphic, monomial, and strictly holomorphic parts, respectively.

In this paper we generalize this situation to coefficients from certain
commutative topological rings, \textit{and} we give explicit formulas for
the decomposition in terms of $\{a_n\}_{n\in\zz}$ and $\{b_n\}_{n\in\zz}$.

The author would like to thank for the support of the Alfr\'ed R\'enyi Institute of Mathematics,
and, in particular, to P\'eter P\'al P\'alfy.

\section{The statement of theorem}

\begin{point}
We say that the topological ring $\mathfrak A$ is a strong polymetric ring if
\begin{itemize}
\item[(a)] its topology is induced by a family of ``seminorms'' $p:\mathfrak A\rightarrow[0,+\infty)$ such that
$p(0)=0$, $p(-X)=p(X)$, $p(X+Y)\leq p(X)+p(Y)$,

\item[(b)] for each ``seminorm'' $p$ there exists a ``seminorm'' $\tilde p$ such that for any $n\in\nn$
$p(X_1\cdot\ldots\cdot X_n)\leq\tilde p(X_1)\cdot\ldots\cdot\tilde p(X_1)$ holds.
\end{itemize}

In what follows $\mathfrak A$ is assumed to be a commutative, separated, sequentially complete strong polymetric ring
with $1\in\mathfrak A$.
\end{point}
\begin{point}
Let $\mathfrak A[z^{-1},z]$ denote the algebra of Laurent series with rapidly decreasing coefficients
with the usual topology.
(We say that a sequence is rapidly decreasing if it is bounded multiplied by any function
of polynomial growth.
The Laurent series $a(z)$ is rapidly decreasing if
the series $p(a_n)$ is rapidly decreasing for any seminorm $p$.)
Let $G=\mathfrak A[z^{-1},z]^\star$ be the unit group, ie.~the group of multiplicatively invertible
elements with topology induced from the topology of $\mathfrak A[z^{-1},z]$ via the
identity inclusion and the multiplicative inverse jointly.

Let
$a(z)\in \mathfrak A[z^{-1},z]^\star$ such that
\[a(z)=\sum_{n\in\zz}a_nz^n,\qquad\qquad\text{and}\qquad\qquad a(z)^{-1}=\sum_{n\in\zz}b_nz^n.\]

(a) We say that $a(z)$ is strictly holomorphic, if $a_0=1$ and $a_n=0$ for $n<0$.
It is easy to see that $a(z)^{-1}$ is also strictly holomorphic.
In fact, the strictly holomorphic Laurent series form a subgroup
$G^+=\mathfrak A[z^{-1},z]^+$ of the unit group. Similar comment applies to
 strictly antiholomorphic Laurent series forming the  subgroup
$G^-=\mathfrak A[z^{-1},z]^-$.

(b) We say that $a(z)$ is an orthogonal Laurent series if for $n,m\in\zz$, $n\neq m$ the identity $a_na_m=0$ holds.
In this case it is immediate that $\Pi^a_n=a_nb_{-n}$ is an idempotent, and the idempotents $\Pi^a_n$
 give a pairwise orthogonal decomposition  of $1$.
Then, $a_n\Pi^a_m=\delta_{n,m}a_n$ holds, for what we can say that $a(z)$ is subordinated to the
orthogonal decomposition $\{\Pi^a_n\}_{n\in\zz}$.
Conversely, if a general $a(z)$ is invertible and it is subordinated to an orthogonal
decomposition then $a(z)$ is also orthogonal, and it is subordinated only to a unique decomposition.
It is easy to see that $\sum_{n\in\zz}\Pi^a_{n}b_{-n} z^{-n}$ is an inverse to $a(z)$,
hence, from the unicity of the inverse $\Pi^a_{-n}b_m=\delta_{n,m}b_m$ follows.
Consequently, $a(z)^{-1}$ is also orthogonal but it is subordinated to
$\{\Pi^a_{-n}\}_{n\in\zz}=\{\Pi^{a^{-1}}_{n}\}_{n\in\zz}$.
Furthermore, if $a_1(z)$ is orthogonal subordinated to $\{\Pi_n^{a_1}\}_{n\in\zz}$
and  $a_2(z)$ is orthogonal subordinated to $\{\Pi_n^{a_2}\}_{n\in\zz}$
then $a(z)=a_1(z)a_2(z)$ will be orthogonal subordinated to $\{\Pi_n\}_{n\in\zz}$
where $\Pi_n=\sum_{m\in\zz}\Pi_m^{a_1}\Pi_{n-m}^{a_2}$.
From that it is easy to see that the orthogonal Laurent series form a subgroup
$\tilde G=\mathfrak A[z^{-1},z]^\sim$ of the unit group.
Clearly, in the case of indecomposable rings, like $\mathfrak A=\mathbb C$,
the orthogonal Laurent series are exactly the (invertible) monomial ones.
\end{point}
Now, we can state our theorem:
\begin{theorem}\label{th:theo}
In accordance to the notation above there is an inner
direct product decomposition
\[G=G^-\times \tilde G\times G^+\]
and the projections $\bo\pi^-$, $\tilde {\bo\pi}$, $\bo\pi^+$ to the various factors, respectively,
are given by explicit infinite algebraic expressions as presented later.
\end{theorem}
This theorem will be proved in Section \ref{sec:proof}.
The proof is done via simple matrix calculations.
In fact, it is a  byproduct of some computations carried out in
\cite{SS}, or, in a more detailed form, in \cite{L}.

\section{Auxiliaries: Matrices}
\begin{point}
Matrices will be considered as linear combinations of elementary matrices $\mathbf e_{n,m}$.
If $S$ is a set of indices then we let $\mathbf 1_{S}=\sum_{s\in S}\mathbf e_{s,s}$.
We will mainly be interested in the cases $S=\zz$ and $S=\zzh\equiv \zz+\frac12$, but we
also define $S^+=\{s\in S\,:\,s>0\}$,  $S^{-}=\{s\in S\,:\,s<0\}$.
When writing  $\zz\times\zz$ matrices  we draw  lines above and below the $0$th row,
and similarly with columns;
when writing  $\zzh\times\zzh$ matrices  we should draw a line between the $\tfrac12$th
and $(-\tfrac12)$th rows, and similarly with columns.
\end{point}
\begin{point}
If $a(z)$ is a Laurent series then we define the multiplication representation matrix as
\[\mathsf U(a(z))=\sum_{n,m\in\mathbb Z}a_{n-m}\mathbf e_{n,m}=
\left[\begin{array}{ccc|c|ccc}
\ddots&\ddots&\ddots&\ddots&\ddots&\ddots&\ddots\\
\ddots&a_0&a_{-1}&a_{-2}&a_{-3}&a_{-4}&\ddots\\
\ddots&a_1&a_0&a_{-1}&a_{-2}&a_{-3}&\ddots\\\hline
\ddots&a_2&a_1&a_0&a_{-1}&a_{-2}&\ddots\\\hline
\ddots&a_3&a_2&a_1&a_0&a_{-1}&\ddots\\
\ddots&a_4&a_3&a_2&a_1&a_0&\ddots\\
\ddots&\ddots&\ddots&\ddots&\ddots&\ddots&\ddots
\end{array}\right].\]

We also define $\mathsf U_{\zzh}(a(z))$ which is the same thing except
using indices $n,m\in\zzh$.
\end{point}
\begin{point}
We say that a $\zz\times\zz$ matrix is a Toeplitz matrix
if it is of shape
\[\mathbf 1_{\zz^-}\mathsf U(a_{\mathrm n}(z))\mathbf 1_{\zz^-}
+\mathbf 1_{\zz^+}\mathsf U(a_{\mathrm p}(z))\mathbf 1_{\zz^+}+A,\]
where $a_{\mathrm n}(z)$, $a_{\mathrm p}(z)$ are rapidly decreasing Laurent series and $A$ is a rapidly
decreasing $\zz\times\zz$ matrix.  $\zzh\times\zzh$ Toeplitz matrices can be defined similarly.
They can be topologized induced from $a_{\mathrm n}(z)$, $a_{\mathrm p}(z)$ and $A$.
The pair $(a_{\mathrm n}(z),a_{\mathrm p}(z))$ is called the symbol of the Toeplitz matrix.
It yields a homomorphism.

In what follows we will always deal with Toeplitz matrices except their coefficients may be not from the
original $\mathfrak A$ but from a larger algebra.
In what follows let $w$ be a formal variable of Laurent series, ie.~whose
coefficients are supposed to be rapidly deceasing.
Similarly, let $t$ be a variable of formal power series, ie.~an infinitesimal variable whose
coefficients are not necessarily rapidly deceasing.
Nevertheless, \textit{if} they are then we can substitute $t=1$.
\end{point}

\begin{point}
We define
\[\mathsf F^R(t,w)=\mathbf 1_{\zz}-tw\mathbf 1_{\zz^-}\mathsf U(z^{-1})-tw^{-1}\mathbf 1_{\zz^+}\mathsf U(z),\]
\[\mathsf F^{R+}(t,w)=\mathbf 1_{\zz}-tw\mathbf 1_{\zz^-}\mathsf U(z^{-1}),\]
\[\mathsf F^{R-}(t,w)=\mathbf 1_{\zz}-tw^{-1}\mathbf 1_{\zz^+}\mathsf U(z).\]

One can notice that
\[\mathsf F^R(t,w)=\mathsf F^{R+}(t,w)\mathsf F^{R-}(t,w)=\mathsf F^{R-}(t,w)\mathsf F^{R+}(t,w).\]
\end{point}
\begin{point}
Let
\[\mathsf U^R(a(z),t,w)=\mathsf F^R(t,w)\mathsf U(a)\mathsf F^R(t,w)^{-1},\]
\[\mathsf U^+(a(z),t,w)=\mathsf F^{R+}(t,w)\mathsf U(a)\mathsf F^{R+}(t,w)^{-1},\]
\[\mathsf U^-(a(z),t,w)=\mathsf F^{R-}(t,w)\mathsf U(a)\mathsf F^{R-}(t,w)^{-1}.\]

For $n>0$ it yields
\begin{multline}
\mathsf U^{R}(z^n,t,w)=\\=
\left[\begin{array}{cccccc|c|cc}
\ddots&&&&&&&&\\
&1&&&&&&&\\\hline
&&1&tw&\cdots&t^{n-1}w^{n-1}&t^nw^n&&\\\hline
&&-tw^{-1}&1-t^2&\cdots&t^{n-2}(1-t^2)w^{n-2}&t^{n-1}(1-t^2)w^{n-1}&&\\
&&&\ddots&\ddots&\vdots&\vdots&&\\
&&&&\ddots&1-t^2&t(1-t^2)w&&\\
&&&&&-tw^{-1}&1-t^2&&\\
&&&&&&&1&\\
&&&&&&&&\ddots
\end{array}\right],\notag\end{multline}
\begin{multline}
\mathsf U^{R}(z^{-n},t,w)=\\=
\left[\begin{array}{cc|c|cccccc}
\ddots&&&&&&&&\\
&1&&&&&&&\\
&&1-t^2&-tw&&&&&\\
&&t(1-t^2)w^{-1}&1-t^2&\ddots&&&&\\
&&\vdots&\vdots&\ddots&\ddots&&&\\
&&t^{n-1}(1-t^2)w^{-n+1}&t^{n-2}(1-t^2)w^{-n+2}&\cdots&1-t^2&-tw&&\\\hline
&&t^nw^{-n}&t^{n-1}w^{-n+1}&\cdots&tw^{-1}&1&&\\\hline
&&&&&&&1&\\
&&&&&&&&\ddots\\
\end{array}\right];\notag
\end{multline}
and
\[\mathsf U^{+}(z^n,t,w)=
\left[\begin{array}{ccccc|c|c}
\ddots&&&&&\\\hline
&1&tw&\cdots &t^{n-1}w^{n-1}&t^nw^n&\\\hline
&&1&\ddots&\ddots&t^{n-1}w^{n-1}&\\
&&&\ddots&\ddots&\vdots&\\
&&&&1&tw&\\
&&&&&1&\\
&&&&&&\ddots
\end{array}\right],\]
\[\mathsf U^{+}(z^{-n},t,w)=
\left[\begin{array}{c|c|ccccc}
\ddots&&&&&&\\
&1&-tw&&&&\\
&&1&\ddots&&&\\
&&&\ddots&\ddots&&\\
&&&&1&-tw&\\\hline
&&&&&1&\\\hline
&&&&&&\ddots
\end{array}\right];\]
and similar matrices show up for $\mathsf U^{-}(z^n,t,w)$ and  $\mathsf U^{-}(z^{-n},t,w)$.
We can see that $\mathsf U^{X}(a(z),t,w)$ is a rapidly decreasing perturbation of $\mathsf U(a(z))$.
Furthermore, we may also notice that it is perfectly legal to substitute $t=1$ into $\mathsf U^{X}(a(z),t,w)$
even if it is illegal to substitute $t=1$ into $\mathsf F^{X}_{\zz}(t,w)^{-1}$.
\end{point}
\begin{point}
We also define
\[\widetilde{\mathsf U}(a(z),w)=
\left[\begin{array}{ccc|c|ccc}
\ddots&\ddots&\vdots&&\vdots&\ddots&\ddots\\
\ddots&a_0&a_{-1}&&-wa_{-2}&-wa_{-3}&\ddots\\
\cdots&a_1&a_0&&-wa_{-1}&-wa_{-2}&\cdots\\\hline
&&&1&&&\\\hline
\cdots&-w^{-1}a_2&-w^{-1}a_1&&a_0&a_{-1}&\cdots\\
\ddots&-w^{-1}a_3&-w^{-1}a_2&&a_1&a_0&\ddots\\
\ddots&\ddots&\vdots&&\vdots&\ddots&\ddots
\end{array}\right].\]
This is also a rapidly decreasing perturbation of $\mathsf U(a(z))$.
\end{point}
\begin{lemma}\label{lem:extract}
We have
\[\mathsf U^R(a(z),1,w)=\widetilde{\mathsf U}(a(z),w)
\left[\begin{array}{ccc|c|ccc}
\ddots&&&&&&\\
&1&&&&&\\
&&1&&&&\\\hline
\cdots&c_{2}(w)&c_{1}(w)&a(w)&c_{-1}(w)&c_{-2}(w)&\cdots\\\hline
&&&&1&&\\
&&&&&1&\\
&&&&&&\ddots
\end{array}\right],\]
where the $c_\lambda(w)$ are rapidly decreasing in $\lambda$ (in both directions).
\begin{proof}
By inspection we see that $\mathsf U^R(a(z),1,w)$
differs from $\widetilde{\mathsf U}(a(z),w)$ only in the $0$th row, and that row can be
brought out by a multiplier as above.
\end{proof}
\end{lemma}

\section{Auxiliaries: Determinants}
\begin{point}
If $A$ is a rapidly decreasing $S\times S$ matrix then $\det(\mathbf 1_S+A)$ can be taken
as usual, as an alternating sum of products associated to finite permutations of $S$.
This yields a  multiplicative operation. In such computations what we have to worry about is that the
sum of the seminorms of the components of
the permanent of the matrix should be convergent,
that provides convergence.
\textit{Sometimes} this naive definition of a determinant  gives a convergent result
even if the argumentum $W$ in not so nice.
Being careful, we use the separate notation $\widetilde\det\, W$ in those cases.

Let us use the notation $\bo[T_1,T_2\bo]=T_1T_2T_1^{-1}T_2^{-1}$ if it makes sense.
If $T_1$ are $T_2$ are invertible Toeplitz matrices (meaning that the inverses are also Toeplitz matrices)
with commuting symbols (as it is obvious in our case) then we see that $\bo[T_1,T_2\bo]$ is a
rapidly decreasing perturbation of the identity matrix.
\end{point}
\begin{lemma}\label{lem:deter}
(a) If $C$ is a rapidly decreasing perturbation of $\mathbf 1_{\zz}$ or $\mathbf 1_{\zzh}$ and $T$ is an
invertible Toeplitz matrix then
\[\det C=\det TCT^{-1}.\]

(b) If $T_1$, $T_2$, $T_3$ are invertible Toeplitz matrices and
$T_1$ commutes with $T_2$ then
\[\det \bo[T_1T_2,T_3\bo]=\det \bo[T_1,T_3\bo]\det \bo[T_2,T_3\bo]. \]

(c) If $T_1$, $T_2$ are invertible Toeplitz matrices, and $T_1'$, $T_2'$ are rapidly decreasing
perturbations, respectively, then
\[\det\left( (T_1'T_2')(T_1T_2)^{-1}\right)=\det T_1'T_1^{-1}\det T_2'T_2^{-1}.\]
\begin{proof}

(a)  If $A$ is a rapidly decreasing matrix and $T$ is any matrix of Toeplitz type then
$\det(\mathbf 1_{S}+AT)=\det(\mathbf 1_{S}+TA)$ is easy to see.
If $T$ is invertible then after putting $AT^{-1}$ to the place of $A$
we see that $\det(\mathbf 1_{S}+A)=\det(\mathbf 1_{S}+TAT^{-1})=\det T(\mathbf 1_{S}+A)T^{-1}$.
(b) Regarding the first part: According to the multiplicativity of the determinant and part (a) we see that
$(\mathrm{LHS})=\det T_2\bo[T_1,T_3\bo]T_2^{-1}\det \bo[T_2,T_3\bo]=(\mathrm{RHS})$.
(c) It follows from  multiplicativity and (a).
\end{proof}
\end{lemma}

\begin{point}
The orthogonal Laurent series $a(z)$ is called orthonormal if $a(1)=1$.
From the discussion earlier it is clear that any  orthogonal $a(z)$ allows
a unique decomposition to a constant and an orthonormal part:
\[a(z)=a(1)\sum_{n\in\zz}\Pi^a_n z^n.\]
There is a very natural way in which orthonormal Laurent series occur:

Suppose that $P$ is a $\zzh\times\zzh$ matrix such that it is a rapidly decreasing perturbation of
$\mathbf 1_{\zzh^-}$ and $P^2=P$. Then, we claim,
\[N_P(z)=\det\left( (\mathbf 1_{\zzh}-P+zP)(\mathbf 1_{\zzh^+}+z\mathbf 1_{\zzh^-})^{-1}  \right)\]
is an orthonormal Laurent series.
This immediately follows from the  the identities $N_P(1)=1$ and $N_P(z_1z_2)=N_P(z_1)N_P(z_2)$;
where $z_1$ and $z_2$ are independent Laurent series variables.
Of these only the second is nontrivial although
it follows from Lemma \ref{lem:deter}.d with the the choices
$T'_i=\mathbf 1_{\zzh}-P+z_iP$ and $T_i=\mathbf 1_{\zzh^+}+z_i\mathbf 1_{\zzh^-}$.
(Conversely, every orthonormal Laurent series occurs in this special form.
Indeed, if $\Pi(z)$ is orthonormal then
$P=\sum^{\mathrm{weak}}_{n\in \zz}\Pi_n\mathbf 1_{\zzh^-+n}$
can be considered. Here the sum was taken in weak sense, ie.~entry-wise in the matrix.
Then one can check that $\Pi(z)=N_P(z)$.)
\end{point}

\section{The proof of the theorem}
\label{sec:proof}
\begin{point}
Let
\[\tilde\pi(a(z),w)=\det \bigl(\mathsf U(a(z))\widetilde{\mathsf U}(a(z),w)^{-1}\bigr),\]
\[\pi^+(a(z),t,w)=\det \left(\mathsf U^+(a(z),t,w)\mathsf U(a(z))^{-1}\right),\]
\[\pi^-(a(z),t,w)=\det \left(\mathsf U^-(a(z),t,w)\mathsf U(a(z))^{-1}\right),\]
\[\pi^\pm(a(z),t,w)=\det \left(\mathsf U^R(a(z),t,w)\mathsf U(a(z))^{-1}\right).\]
\end{point}
\begin{lemma}\label{lem:X}The following identities hold:
\begin{itemize}
\item[(a)]
\[\pi^\pm(a(z),t,w)=\pi^+(a(z),t,w)\pi^-(a(z),t,w).\]
\item[(b)]
\[\pi^\pm(a(z),1,w)=\pi^+(a(z),1,w)\pi^-(a(z),1,w).\]
\item[(c)]
\[\pi^\pm(a(z),1,w)=a(w)\tilde \pi(a(z),w)^{-1}.\]
\item[(d)]
\begin{multline}
\qquad\qquad\tilde\pi(a(z),w)\tilde\pi(a(z),1)^{-1}=
\det(\widetilde{\mathsf U}(a(z),1) \widetilde{\mathsf U}(a(z),w)^{-1})=\\
=\det\left(\mathsf U_{\zzh}(a(z))(w\mathbf 1_{\zzh^-}+\mathbf 1_{\zzh^+})
\mathsf U_{\zzh}(a(z)^{-1})(w\mathbf 1_{\zzh^-}+\mathbf 1_{\zzh^+})^{-1}\right).\qquad\qquad\notag\end{multline}
\end{itemize}
\begin{proof}
(a) It follows from
Lemma \ref{lem:deter}.b with the choices $T_1=\mathsf F^{R+}(t,w)$, $T_2=\mathsf F^{R-}(t,w)$, $T_3=\mathsf U(a(z))$.
(b) We can substitute $t=1$ into (a).
(c) It follows from Lemma \ref{lem:extract}.
(d) The first equality follows immediately from the multiplicativity of the determinant.
Regarding the second, take $\widetilde{\mathsf U}(a(z),1) \widetilde{\mathsf U}(a(z),w)^{-1}$.
Remove the single element $1$ in the center. Then we obtain a $(\zz\setminus\{0\})\times (\zz\setminus\{0\})$
matrix.
Relabel the elementary matrices $\mathbf e_{n,m}$ into $\mathbf e_{n-\frac12\sgn n,m-\frac12\sgn m}$.
Conjugate by $\mathbf 1_{\zzh^-}-\mathbf 1_{\zzh^+}$. Then the determinant does not change, yet
we obtain $ \mathsf U_{\zzh}(a(z))(w\mathbf 1_{\zzh^-}+\mathbf 1_{\zzh^-})
\mathsf U_{\zzh}(a(z)^{-1})(w\mathbf 1_{\zzh^-}+\mathbf 1_{\zzh^-})^{-1}$
at the end.
\end{proof}
\end{lemma}
\begin{lemma}The following identities/statements hold:
\begin{itemize}
\item[(a)] \[a(w)=\pi^-(a(z),1,w)\tilde\pi(a(z),w)\pi^+(a(z),1,w).\]
\item[(b)] \[\pi^-(a(z),1,w)\in\mathfrak A[w,w^{-1}]^-,\]
\[\tilde\pi(a(z),w)\in\mathfrak A[w,w^{-1}]^\sim,\]
\[\pi^+(a(z),1,w)\in\mathfrak A[w,w^{-1}]^+.\]
\item[(c)]
\[\pi^+(a_1(z)a_2(z),1,w)=\pi^+(a_1(z),1,w)\pi^+(a_2(z),1,w),\]
\[\tilde\pi(a_1(z)a_2(z),w)=\tilde \pi(a_1(z),w)\tilde \pi(a_2(z),w),\]
\[\pi^-(a_1(z)a_2(z),1,w)=\pi^-(a_1(z),1,w)\pi^-(a_2(z),1,w).\]
\item[(d)]
\[\tilde\pi(a(z),w)=\pi^+(a(z),1,w)=1\quad\text{if}\quad a(z)\in\mathfrak A[z,z^{-1}]^-,\]
\[\pi^-(a(z),1,w)=\pi^+(a(z),1,w)=1\quad\text{if}\quad a(z)\in\mathfrak A[z,z^{-1}]^\sim,\]
\[\pi^-(a(z),1,w)=\tilde\pi(a(z),w)=1\quad\text{if}\quad a(z)\in\mathfrak A[z,z^{-1}]^+.\]
\end{itemize}
\begin{proof}

(a) It is a combination of Lemma \ref{lem:X}.b and c.
(b): The first and third lines follow from the special form of the matrices involved,
which can be understood as special power series in $w$ and $w^{-1}$ respectively.
The second line follows from Lemma \ref{lem:X}.d and the discussion about orthonormal Laurent series
applied with the choice $P=\mathsf U(a(z)^{-1})\mathbf 1_{\zzh^-}\mathsf U(a(z))$.
(c) They follow from Lemma \ref{lem:deter}.c.
(d) They follow from that we take the determinant of
upper or lower triangular matrices with 1's in the diagonals.
\end{proof}
\end{lemma}

\begin{point}
Then we can define the operations $\bo\pi^-$, $\tilde{\bo\pi}$, $\bo\pi^+$, respectively, for $a(z)$
by taking $\pi^-(a(z),1,w)$, $\tilde\pi(a(z),w)$, $\pi^+(a(z),1,w)$, respectively,
and replacing $w$ by $z$ in the result.
Now, the content of the previous lemma is exactly that these operations
serve an inner direct product decomposition of as indicated in Theorem \ref{th:theo}.
By this we have proved the theorem.
\end{point}

\section{Explicit formulas}
\begin{point}
The operations $\bo\pi^-$, $\tilde{\bo\pi}$, $\bo\pi^+$, or rather
$\pi^-(a(z),1,w)$, $\tilde\pi(a(z),w)$, $\pi^+(a(z),1,w)$
are already given in a sufficiently explicit form.
Yet, it is useful to have transparent and uniform formulas as much as it is possible.
\end{point}
\begin{lemma}\label{lem:exp}
Suppose that $A$ is a rapidly decreasing $\zz\times\zz$ matrix. Then
\[\widetilde{\det}\left( \mathsf F^{R+}(t,w)+A\right)=
\det\left( (\mathsf F^{R+}(t,w)+A)\mathsf F^{R+}(t,w)^{-1}\right).\]
Similar statement holds with $\mathsf F^{R-}(t,w)$.
\newpage
\begin{proof}
For finite matrices $A$ the statement is relatively easy to see.
Due to the special form of the matrices the permanents  are sufficiently
controlled, hence the statement extends.
\end{proof}
\end{lemma}
\begin{prop}
\begin{multline}
\pi^+(a(z),1,w)=\widetilde{\det}\left(\mathsf U(a(z)^{-1})\mathsf F^{R+}(1,w)\mathsf U(a(z))\right)=\\
=\widetilde{\det}\left(\mathbf 1_{\zz}-w\mathsf U(a(z)^{-1})\mathbf 1_{\zz^-}\mathsf U(a(z))\mathsf U(z^{-1})\right),
\notag\end{multline}\vspace{-6mm}
\begin{multline}
 \tilde\pi(a(z),w)=\det\left(\mathsf U(a(z)\widetilde{\mathsf U}(a(z),1)^{-1}\right)\cdot\\\cdot\det\left(
\mathsf U_{\zzh}(a(z)^{-1})(w\mathbf 1_{\zzh^-}+\mathbf 1_{\zzh^-})^{-1}
\mathsf U_{\zzh}(a(z))(w\mathbf 1_{\zzh^-}+\mathbf 1_{\zzh^-})
\right),\notag\end{multline}\vspace{-6mm}
\begin{multline}
\pi^-(a(z),1,w)=\widetilde{\det}\left(
\mathsf U(a(z)^{-1})\mathsf F^{R-}(1,w)\mathsf U(a(z))
\right)=\\=\widetilde{\det}\left(
\mathbf 1_{\zz}-w^{-1}\mathsf U(a(z)^{-1})\mathbf 1_{\zz^+}\mathsf U(a(z))\mathsf U(z)\right) .
\notag\end{multline}
\begin{proof}
The second equality follows from Lemma \ref{lem:X}.d, while the other ones
follow from Lemma \ref{lem:exp}.
\end{proof}
\end{prop}
\begin{point}
The point is that the matrices
$\mathsf U(a(z)^{-1})\mathbf 1_{\zz^-}\mathsf U(a(z))$,
$\mathsf U(a(z)^{-1})\mathbf 1_{\zz^+}\mathsf U(a(z))$
are relatively transparent perturbations of $\mathbf 1_{\zz^-}$ and $\mathbf 1_{\zz^+}$.
In particular, if $a(z)$ is finite then
these perturbations restrict to finitely many columns, hence, as a result of the block triangular form,
effectively reducing the computations to finite matrices.
Cf.~like when
\begin{multline}
\mathsf U(a(z)^{-1})\mathsf F^{R+}(1,w)\mathsf U(a(z))
=\mathbf 1_{\zz}-w\mathsf U(a(z)^{-1})\mathbf 1_{\zz^-}\mathsf U(a(z))\mathsf U(z^{-1})=\\
=\left[\begin{array}{ccc|ccc|ccc}
\ddots&\ddots&&\vdots&\vdots&\vdots&&&\\
&1&-w&*&*&*&&&\\
&&1&*&*&*&&&\\\hline
&&&*&*&*&&&\\
&&&*&*&*&&&\\
&&&*&*&*&&&\\\hline
&&&*&*&*&1&&\\
&&&*&*&*&&1&\\
&&&\vdots&\vdots&\vdots&&&\ddots
\end{array}\right];
\notag\end{multline}
and taking determinant reduces to the block in the middle.
Even there, the determinant will have the shape of a characteristic polynomial
$\det(\mathbf 1_S-wB)$.
\end{point}

\end{document}